\documentclass[12pt]{article}  
\usepackage{amsthm,amsmath,amsfonts,epsfig,amssymb,latexsym}
\usepackage{mathrsfs}
\usepackage[T1]{fontenc}
\usepackage{graphicx}
\usepackage{enumerate}
\usepackage{xy,xypic}
\usepackage{graphics}
\usepackage{footnote}
\usepackage{scalefnt}
\usepackage{tikz}
\usetikzlibrary{shapes,arrows,fit,calc,positioning}
\usetikzlibrary{matrix}
\xyoption{all}
\title{{\Large Some perspective on Homotopy obstructions} } 
\author{
 Satya Mandal\footnote{Partially supported by a General Research Grant (no 2301857) from U. of Kansas}
 ~~~and ~~Bibekananda Mishra
 \\ 
{\small University of Kansas, Lawrence, Kansas 66045, USA}\\
{\small {\it  mandal@ku.edu, bibekanandamishra@ku.edu} 
  }\\
 } 
\begin{document}
\renewcommand{\baselinestretch}{1.255}
\setlength{\parskip}{1ex plus0.5ex}
\date{11 November 2018}

\newtheorem{theorem}{Theorem}[section]
\newtheorem{proposition}[theorem]{Proposition}
\newtheorem{lemma}[theorem]{Lemma}
\newtheorem{corollary}[theorem]{Corollary}
\newtheorem{construction}[theorem]{Construction}
\newtheorem{notations}[theorem]{Notations}
\newtheorem{question}[theorem]{Question}
\newtheorem{example}[theorem]{Example}
\newtheorem{definition}[theorem]{Definition} 
\newtheorem{conjecture}[theorem]{Conjecture} 
\newtheorem{remark}[theorem]{Remark} 
\newtheorem{statement}[theorem]{Statement}

\newcommand{\iso}{\stackrel{\sim}{\longrightarrow}}

\newcommand{\sur}{\twoheadrightarrow}
\newcommand{\bD}{\begin{definition}}
\newcommand{\eD}{\end{definition}}
\newcommand{\bP}{\begin{proposition}}
\newcommand{\eP}{\end{proposition}}
\newcommand{\bL}{\begin{lemma}}
\newcommand{\eL}{\end{lemma}}
\newcommand{\bT}{\begin{theorem}}
\newcommand{\eT}{\end{theorem}}
\newcommand{\bC}{\begin{corollary}}
\newcommand{\eC}{\end{corollary}} 
\newcommand{\eop}{\hfill \rule{2mm}{2mm}}
\newcommand{\pf}{\noindent{\bf Proof.~}}
\newcommand{\PD}{\text{proj} \dim}
\newcommand{\lra}{\longrightarrow}
\newcommand{\hra}{\hookrightarrow}
\newcommand{\llra}{\longleftrightarrow}
\newcommand{\Lra}{\Longrightarrow}
\newcommand{\Llra}{\Longleftrightarrow}
\newcommand{\bE}{\begin{enumerate}}
\newcommand{\eE}{\end{enumerate}}
\newcommand{\Sets}{\underline{{\mathrm Sets}}}
\newcommand{\Sch}{\underline{{\mathrm Sch}}}
\newcommand{\ForMe}{\noindent\TCP{{\bf Remarks To Be Removed:~}}}
\newcommand{\pic}{The proof is complete.}
\newcommand{\tcp}{This completes the proof.}

\def\spec#1{\mathrm{Spec}\left(#1\right)}
\def\m{\mathfrak {m}}
\def\CA{\mathcal {A}}
\def\CB{\mathcal {B}}
\def\CP{\mathcal {P}}
\def\CC{\mathcal {C}}
\def\CD{\mathcal {D}}
\def\CE{\mathcal {E}}
\def\CF{\mathcal {F}}
\def\CE{\mathcal {E}}
\def\CG{\mathcal {G}}
\def\CH{\mathcal {H}}
\def\CI{\mathcal {I}}
\def\CJ{\mathcal {J}}
\def\CK{\mathcal {K}}
\def\CL{\mathcal {L}}
\def\CM{\mathcal {M}}
\def\CN{\mathcal {N}}
\def\CO{\mathcal {O}}
\def\CP{\mathcal {P}}
\def\CQ{\mathcal {Q}}
\def\CR{\mathcal {R}}
\def\CS{\mathcal {S}}
\def\CT{\mathcal {T}}
\def\CU{\mathcal {U}}
\def\CV{\mathcal {V}}
\def\CW{\mathcal {W}}
\def\CX{\mathcal {X}}
\def\CY{\mathcal {Y}}
\def\CZ{\mathcal {Z}}

\newcommand{\smallcirc}[1]{\scalebox{#1}{$\circ$}}
\def\BA{\mathbb {A}}
\def\BB{\mathbb {B}}
\def\BC{\mathbb {C}}
\def\BD{\mathbb {D}}
\def\BE{\mathbb {E}}
\def\BF{\mathbb {F}}
\def\BG{\mathbb {G}}
\def\BH{\mathbb {H}}
\def\BI{\mathbb {I}}
\def\BJ{\mathbb {J}}
\def\BK{\mathbb {K}}
\def\BL{\mathbb {L}}
\def\BM{\mathbb {M}}
\def\BN{\mathbb {N}}
\def\BO{\mathbb {O}}
\def\BP{\mathbb {P}}
\def\BQ{\mathbb {Q}}
\def\BR{\mathbb {R}}
\def\BS{\mathbb {S}}
\def\BT{\mathbb {T}}
\def\BU{\mathbb {U}}
\def\BV{\mathbb {V}}
\def\BW{\mathbb {W}}
\def\BX{\mathbb {X}}
\def\BY{\mathbb {Y}}
\def\BZ{\mathbb {Z}}

\newcommand{\TCP}{\textcolor{purple}}
\newcommand{\TCM}{\textcolor{magenta}}
\newcommand{\TCR}{\textcolor{red}}
\newcommand{\TCB}{\textcolor{blue}}
\newcommand{\TCG}{\textcolor{green}}

\def\SA{\mathscr {A}}
\def\SB{\mathscr {B}}
\def\SC{\mathscr {C}}
\def\SD{\mathscr {D}}
\def\SE{\mathscr {E}}
\def\SF{\mathscr {F}}
\def\SG{\mathscr {G}}
\def\SH{\mathscr {H}}
\def\SI{\mathscr {I}}
\def\SJ{\mathscr {J}}
\def\SK{\mathscr {K}}
\def\SL{\mathscr {L}}
\def\SN{\mathscr {N}}
\def\SO{\mathscr {O}}
\def\SP{\mathscr {P}}
\def\SQ{\mathscr {Q}}
\def\SR{\mathscr {R}}
\def\SS{\mathscr {S}}
\def\ST{\mathscr {T}}
\def\SU{\mathscr {U}}
\def\SV{\mathscr {V}}
\def\SW{\mathscr {W}}
\def\SX{\mathscr {X}}
\def\SY{\mathscr {Y}}
\def\SZ{\mathscr {Z}}

\def\bfA{{\bf A}}
\def\bfB{{\bf B}} 
\def\bfC{{\bf C}} 
\def\bfD{{\bf D}} 
\def\bfE{{\bf E}} 
\def\bfF{{\bf F}} 
\def\bfG{{\bf G}} 
\def\bfH{{\bf H}} 
\def\bfI{{\bf I}} 
\def\bfJ{{\bf J}} 
\def\bfK{{\bf K}} 
\def\bfL{{\bf L}} 
\def\bfM{{\bf M}} 
\def\bfN{{\bf N}} 
\def\bfO{{\bf O}} 
\def\bfP{{\bf P}} 
\def\bfQ{{\bf Q}} 
\def\bfR{{\bf R}} 
\def\bfS{{\bf S}} 
\def\bfT{{\bf T}} 
\def\bfU{{\bf U}} 
\def\bfV{{\bf V}} 
\def\bfW{{\bf W}} 
\def\bfX{{\bf X}} 
\def\bfY{{\bf Y}} 
\def\bfZ{{\bf Z}} 

\maketitle

\section{Introduction} 
Throughout this article $A$ will denote a commutative noetherian ring, with $\dim A=d\geq 2$, and
$A[T]$ will denote the polynomial
ring in one varaible $T$. Also,
$P$ will denote a projective $A$-module with $rank(P)=n$.
In \cite{MM1} we considered the  Homotopy obstruction sets $\pi_0\left({\mathcal LO}(P)\right)$. In this article, we 
provide some perspective on these sets $\pi_0\left({\mathcal LO}(P)\right)$, by proving that some of them are isomorphic and by defining  natural 
set theoretic maps to Chow groups.

Assume  $A$ is a regular ring containing a field $k$, with $1/2\in k$.
In \cite{MM1}, it was established that,
 if $2n\geq d+2$, then $\pi_0\left({\mathcal LO}(P)\right)$ has a natural structure of a 
abelian monoid.
Further, if $P\cong P_0\oplus A$, then $\pi_0\left({\mathcal LO}(P)\right)$ is an abelian group. 
Let $P$, $Q$ be two projective $A$-module with $rank(P)=rank(Q)=d$ and let $\iota: \Lambda^d Q\iso  \Lambda^d P$ be an isomorphism of the determinants.
Under the same regularity hypotheses, in this article, we prove that there is a natural monoid
 isomorphism $\chi(\iota): \pi_0\left({\mathcal LO}(P)\right) \iso \pi_0\left({\mathcal LO}(Q)\right)$. Therefore, when $rank(P)=d$, we have the following:
 
   \bE
   \item There is a natural monoid isomorphism 
 $$
  \pi_0\left({\mathcal LO}(P)\right)\cong \pi_0\left({\mathcal LO}(\Lambda^d P\oplus A^{d-1})\right)
  $$

  \item Since $\pi_0\left({\mathcal LO}(\Lambda^d P\oplus A^{d-1})\right)$ is an abelian group, so is $\pi_0\left({\mathcal LO}(P)\right)$,
  for all such $P$.
    \item We conclude that, for $\pi_0\left({\mathcal LO}(P)\right)$ to have a group structure, it is not necessary that $P\cong P_0\oplus A$. 
  \item If $A$ is essentially smooth over an infinite perfect field $k$, with $d\geq 3$, 
 it follows that $ \pi_0\left({\mathcal LO}(P)\right)\cong  E^d\left(A, \Lambda^d P\right)$
 where $E^d\left(A, \Lambda^d P\right)$ denotes the Euler class group, as defined in \cite{BS1, BS2, MY, MM1}.
\eE

The above were discussed in Section \ref{mainCompa}. 
In Section \ref{ChowSec}, we compare $\pi_0\left({\mathcal LO}(P)\right)$ with the Chow groups, as follows.
 Let $A$ be a regular  ring, over an infinite field $k$, with $1/2\in k$ and $\dim A=d$. Let $P$ be any projective $A$-module, with $rank(P)=n$. 
In Section \S \ref{ChowSec}, we  establish a natural set theoretic map $\pi_0\left({\mathcal LO}(P)\right) \lra CH^n(A)$, where $CH^n(A)$ denotes the 
Chow group of codimension $n$ cycles \cite{F}. In section \ref{HomoDesc}, we record and formalize an expected alternate  description of 
$\pi_0\left({\mathcal LO}(P)\right)$, which was instrumental in establishing the results in other sections.

For all unexplained notations and definitions, the readers are referred \cite{MM1}.
 In particular, for an $A$-module $M$, denote
$M[T]:=M\otimes A[T]$. Likewise, for a homomorphism $f:M\lra N$ of  $A$-modules, $f[T]:=f\otimes A[T]$.


\section{Alternate Description of  Obstructions $\pi_0\left({\mathcal LO}(P) \right)$} \label{HomoDesc}
Let $A$ be a noetherian commutative  ring, with $\dim A=d$ and $P$ be a projective $A$-module with $rank(P)=n$. In \cite{MM1}, we defined 
Nori Homotopy  Obstruction sets $\pi_0\left({\mathcal LO}(P)\right)$, and gave several other description of the same. In this section, 
we recall the essential part of the definition of $\pi_0\left({\mathcal LO}(P)\right)$, and give one more description of the same, as follows.

\bD\label{defpizero2Oct}{\rm 
Let $A$ be a noetherian commutative  ring, with $\dim A=d$ and $P$ be a projective $A$-module with $rank(P)=n$.
 By a {\bf local $P$-orientation}, we mean 
a pair $(I, \omega)$ where $I$ is an ideal of $A$ and $\omega:P \sur \frac{I}{I^2}$ is a surjective homomorphism.
We will use the same notation $\omega$ for the map $\frac{P}{IP} \sur \frac{I}{I^2}$, induced by $\omega$. 
 A local {\bf local $P$-orientation} will simply be 
 referred to as a {\bf local orientation}, when $P$ is understood.
%
Denote 
\begin{equation}\label{4LOetc}
\left\{
\begin{array}{l}
{\mathcal LO}(P)=\left\{(I, \omega): (I, \omega)~{\rm is~a~local}~P~{\rm orientation} \right\}\\
{\mathcal LO}^{\geq n}(P)=\left\{(I, \omega)\in {\mathcal LO}(P): height(I)\geq n  \right\}\\
\end{array}
\right.
\end{equation}
%
 In \cite{MM1}, the obstruction set $\pi_0\left({\mathcal LO}(P)\right)$ was defined, by the push forward diagram:
\begin{equation}\label{pizeroLOpush}
\diagram
{\mathcal LO}(P[T])\ar[r]^{T=0}\ar[d]_{T=1} &{\mathcal LO}(P)\ar[d]\\
{\mathcal LO}(P)\ar[r] & \pi_0\left({\mathcal LO}(P)\right)\\ 
\enddiagram, 
%
{\rm in}\quad \underline{Sets}.
\end{equation}
While one would like to define $\pi_0\left({\mathcal LO}^{\geq n}(P)\right)$ similarly, note that substitution $T=0, 1$ would not yield any map 
from   ${\mathcal LO}^{\geq n}(P[T])$ to ${\mathcal LO}^{\geq n}(P)$. However, note that the definition of $\pi_0\left({\mathcal LO}(P)\right)$ 
by push forward diagram (\ref{pizeroLOpush}), is only an alternate way of saying the following:
\bE
\item 
%
%
For $(I_0, \omega_0), (I_0, \omega_0)\in {\mathcal LO}(P)$, 
we write 
$(I_0, \omega_0)\sim (I_0, \omega_0)$, if there is an $(I, \omega)\in {\mathcal LO}(P[T])$, if $(I(0), \omega(0))=(I_0, \omega_0)$ and 
$(I(1), \omega(1))=(I_1, \omega_1)$. 
Now $\sim$ generates a chain equivalence relation on ${\mathcal LO}(P)$, which we call the chain homotopy relation.
\item The above definition (\ref{pizeroLOpush}) means,
$\pi_0\left({\mathcal LO}(P)\right)$  is the set of all equivalence classes in  ${\mathcal LO}(P)$.
\eE

The restriction of the  relation $\sim$ on  ${\mathcal LO}^{\geq n}(P) \subseteq {\mathcal LO}(P)$,
generates a chain homotopy relation on  ${\mathcal LO}^{\geq n}(P)$. 
Define $\pi_0\left({\mathcal LO}^{\geq n}(P)\right)$ to be the set of all equivalence classes in ${\mathcal LO}^{\geq n}(P)$.
It follows, that there is natural map
$$
\varphi: \pi_0\left({\mathcal LO}^{\geq n}(P)\right) \lra \pi_0\left({\mathcal LO}(P)\right)
$$
}
\eD
%
\bP\label{transIso}
Let $A$ and $P$ be as in (\ref{defpizero2Oct}). Then,
\bE
\item 
The map $\varphi$ is surjective.
\item \label{bijRegHyp} If $A$ is a regular ring  containing is field $k$, with $1/2\in k$, then $\varphi$ a bijection.
\item\label{teenRee} Let $A$ be as in (\ref{bijRegHyp}). Then $\sim$ is  an equivalence relation on $\pi_0\left({\mathcal LO}^{\geq n}(P)\right)$.

\eE
 
\eP
\pf Let $x=[(I, \omega)] \in  \pi_0\left({\mathcal LO}(P)\right)$. By application of the Involution operation \cite[Section 5]{MM1} twice, we can assume $height(I)\geq n$.
So, $(I, \omega) \in {\mathcal LO}^{\geq n}(P)$. This establishes that $\varphi$ is surjective.

Now, assume $A$ is as in (\ref{bijRegHyp}).  
Suppose  $x_0, x_1\in \pi_0\left({\mathcal LO}^{\geq n}(P)\right)$, and $\varphi(x_0) =\varphi(x_1)$. Then, form $i=0, 1$, we have
$x_i=[(I_i, \omega_i)]$,
 for some $(I_i, \omega_i) \in {\mathcal LO}^{\geq n}(P)$. By \cite[Corollary 3.2]{MM1}, there is a homotopy 
$H(T)=(J, \Omega)\in {\mathcal LO}(P[T])$ such that $H(0)=(I_0, \omega_0)$ and $H(1)=(I_1, \omega_1)$.  Therefore, $x_0=x_1$. 
This establishes (\ref{bijRegHyp}) and (\ref{teenRee}) follows by the same argument.
\pic $\eop$

%
\begin{remark}\label{moviongRem}{\rm 
Use the notations as in Definition \ref{defpizero2Oct}. It would be worthwhile recording the following,
\bE
\item Suppose $(I_0, \omega_0), (I_1, \omega_1)\in {\mathcal LO}^{\geq n}(P)$ and $(I_0, \omega_0)\sim (I_1, \omega_1)$. By defintion there is homotopy
$H(T)=(I, \omega)\in {\mathcal LO}(P[T])$ such that $H(0)= (I_0, \omega_0)$ and $H(1)=(I_1, \omega_1)$. Note,
by moving Lemma argument, similar to \cite[Lemma 4.5]{MM1}, we can assume that $H(T)\in {\mathcal LO}^{\geq n}(P[T])$.
%
\item Assume $A$ is a Cohen Macaulay ring. Then, ${\mathcal LO}^{\geq n}(P)$ is in bijection with the set
$$
\left\{(I, \omega) \in {\mathcal LO}(P): height(I)=n,~\omega: \frac{P}{IP} \iso \frac{I}{I^2}~{\rm is~an~isomorphism}\right\}
\cup \left\{(A, 0)\right\}
$$
 
\eE
}
\end{remark}

\section{Comparison of $\pi_0\left({\mathcal LO}(P)\right)$: top rank case}\label{mainCompa}
In this section, under usual regularity hypotheses, 
we prove that,  the Homotopy obstruction sets $\pi_0\left({\mathcal LO}(P)\right)$ are isomorphic, when $rank(P)=d$

\bD\label{defnatMap}{\rm
Suppose $A$ is a regular ring, containing  field $k$, with $1/2\in k$.
Let $P, Q$ be two projective $A$-modules with $rank(P)=rank(Q)=\dim A=d\geq 2$. 
Assume they have isomorphic determinant and $\iota: \Lambda^d Q\iso \Lambda^dP$ be an isomorphism.
Then, we define a natural map 
$$
\chi(\iota): \pi_0\left({\mathcal LO}(P)\right) \lra \pi_0\left({\mathcal LO}(Q)\right) 
$$
By Proposition \ref{transIso}, there is a bijection $\pi_0\left({\mathcal LO}^{\geq n}(P)\right)\iso \pi_0\left({\mathcal LO}(P)\right)$. So,
we would define a
$$
\chi(\iota): \pi_0\left({\mathcal LO}^{\geq n}(P)\right) \lra \pi_0\left({\mathcal LO}^{\geq n}(Q)\right) 
\qquad \qquad {\rm as ~follows:}
$$
\bE
\item\label{chiZero} Let $(I, \omega)\in {\mathcal LO}^{\geq n}(P)$, with $height(I)=d$. Use the same notation $\omega:\frac{P}{IP}\iso \frac{I}{I^2}$,
for the map (isomorphism), induced by $\omega: P \sur \frac{I}{I^2}$.
%
Let 
$$
\alpha:\frac{Q}{IQ} \iso \frac{P}{IP} \quad {\rm be~ an~ isomorphism,}\quad\ni\quad  \det(\alpha)=\iota\otimes \frac{A}{I}.
$$
Let $\tilde{\omega}(\alpha)$  be defined by the commutative diagram
$$
\diagram 
Q\ar[r]\ar@{-->}[drr]_{\tilde{\omega}(\alpha)} & \frac{Q}{IQ}\ar[r]^{\alpha} & \frac{P}{IP} \ar[d]^{\omega}\\
&& \frac{I}{I^2}\\
\enddiagram
\quad {\rm and ~define}\quad 
\left\{\begin{array}{l}
\chi_0(\iota)(I, \omega)=(I, \tilde{\omega}(\alpha))\in {\mathcal LO}(Q)\\
\chi_1(\iota)(I, \omega)= [\chi_0(I, \omega)]\in \pi_0\left({\mathcal LO}(Q)\right).\\
\end{array}
\right.
$$
\item \label{chiIota}
For $x\in \pi_0\left({\mathcal LO}(P) \right)$, we can write $x=[(I, \omega)]$, for some $(I, \omega)\in {\mathcal LO}(P)$, with
$height(I)\geq d$. We define,
\begin{equation}\label{DefChiIOTA}
\chi(\iota)(x) =
\left\{
\begin{array}{lr}
\chi_1(\iota)(I, \omega)\in \pi_0\left({\mathcal LO}(Q) \right) & if~height(I)=d\\
\left[(A, 0)\right] \in  \pi_0\left({\mathcal LO}(Q) \right) & if~I=A\\
\end{array}
\right.
\end{equation}

\eE
Subsequently, we prove that $\chi_0(\iota)(I, \omega)$ is independent of the choice of $\alpha$, and 
also $\chi(\iota):\pi_0\left({\mathcal LO}(P)\right) \lra \pi_0\left({\mathcal LO}(Q)\right)$ is well defined.
}
\eD
\bL\label{IndOfAlpha}
With notations as in (\ref{defnatMap}),
 let $(I, \omega)\in {\mathcal LO}(P)$, with $height(I)=d$.
Then, $\chi_0(\iota)(I, \omega)$ is independent of the choice of $\alpha$.
\eL
\pf   
Let $\alpha$ be as in (\ref{defnatMap}), and 
$$
\beta:\frac{Q}{IQ} \iso \frac{P}{IP} \quad {\rm be~ another~ isomorphism,}\quad\ni\quad  \det(\beta)=\iota\otimes \frac{A}{I}.
$$
Then, $\det(\alpha^{-1}\beta)=1$. Since $\dim\left(\frac{A}{I}\right)=0$, the map 
 $\alpha^{-1}\beta$ is an elementary matrix, with respect to any choice of basis of $\frac{Q}{IQ}$.
Therefore, there is an elementary automorphism 
 $\gamma(T)\in EL\left(\frac{Q[T]}{IQ[T]}\right)$,
such that $\gamma(0)=Id$ and $\gamma(1)= \alpha^{-1}\beta$. Consider $(I[T], H(T))\in {\mathcal LO}(Q[T])$, where $H(T)$
is defined by the commutative diagram:
$$
\diagram 
Q[T]\ar@{-->}[drr]_{H(T))} \ar[r] & \frac{Q[T]}{IQ[T]} \ar[r]^{\gamma(T)} &  \frac{Q[T]}{IQ[T]} \ar[d]
\ar[r]^{\alpha[T]}& \frac{P[T]}{IP[T]}\ar[dl]^{\omega[T]}\\
&& \frac{I[T]}{I[T]^2}&\\
\enddiagram 
$$
Then with $T=0, 1$, we have
$$
\diagram 
Q\ar@{-->}[drr]_{H(0)} \ar[r] & \frac{Q}{IQ} \ar[r]^{1} &  \frac{Q}{IQ} \ar[d]^{\omega\alpha}\\
&& \frac{I}{I^2}\\
\enddiagram 
\qquad {\rm and}\qquad 
\diagram 
Q\ar@{-->}[drr]_{H(1)} \ar[r] & \frac{Q}{IQ} \ar[r]^{\alpha^{-1}\beta} &  \frac{Q}{IQ} \ar[d]^{\omega\alpha}\\
&& \frac{I}{I^2}\\
\enddiagram 
$$
This completes the proof. $\eop$

\bL\label{DefofChi}{\rm 
Use the notations as in Definition \ref{defnatMap}. Then,  the definition (\ref{DefChiIOTA}) 
of the  map $\chi(\iota):\pi_0\left({\mathcal LO}(P)\right) \lra \pi_0\left({\mathcal LO}(Q)\right)$, 
 is well defined.

}
\eL
\pf Let $x=[(I, \omega)]= [(I_1, \omega_1)] \in  \pi_0\left({\mathcal LO}(P) \right)$, with $height(I)\geq d$ and $height(I_1)\geq d$. 
Since homotopy induces an equivalence relation \cite[Corollary 3.2 ]{MM1}, there is a homotopy ${\CH}(T)=(J, \omega_J) \in {\mathcal LO}(P[T])$ such that 
$$
{\CH}(0)=(J(0), \omega_J(0))=(I, \omega),\qquad {\rm and}\qquad {\CH}(1)=(J(1), \omega_J(1))=(I_1, \omega_1)
$$
By Moving Lemma \cite[Lemma 4.5]{MM1}, we can assume that $height(J)\geq d$ ({\it In deed, $height(J)=n$, unless $I=I_1=A$}).
Since $\dim\left(\frac{A[T]}{J}\right)=1$, it follows, there is an isomorphism 
$$
\alpha: \frac{Q[T]}{JQ[T]} \iso \frac{P[T]}{JP[T]} 
\quad\ni\quad  \det(\alpha)=\iota\otimes \frac{A[T]}{J}
$$
 Now define ${\CH}(T)=(J, \Omega)$, where $\Omega$ is defined as follows:
$$
\diagram 
Q[T] \ar[r] \ar@{-->}[rdr]_{\Omega}& \frac{Q[T]}{JQ[T]} \ar[r]^{\alpha} & \frac{P[T]}{JP[T]} \ar[d]^{\omega_J}\\
&& \frac{J}{J^2}\\
\enddiagram
$$
It follows 
$$
\chi_1(I, \omega)= [(I, \tilde{\omega}(\alpha(0)))] =[(J(0), \Omega(0))]= [(J(1), \Omega(1))]= [(I_1, \tilde{\omega}_1(\alpha(1)))]=\chi_1(I_1, \omega_1).
$$
\pic $\eop$

\bP\label{MonoidHomo}
Use the notations as in Definition \ref{defnatMap}. Then, the map $\chi(\iota)$ is a monoid homomorphism. 
\eP
\pf Let $x=[(I_1, \omega_1)], y=[(I_2, \omega_2)\in \pi_0\left({\mathcal LO}(P)\right)$. We can assume $I_1+I_2=A$ and $height(I_i)\geq d$ for $i=1, 2$.

If $I_1=A$ or $I_2=A$, then by definition (\ref{DefChiIOTA}), we have $\chi(\iota)(x+y)=\chi(\iota)(x)+\chi(\iota)(y)$.
So, assume $height(I_1)=height(I_2)=d$. 
For $i=1, 2$ let
$$
\alpha_i:\frac{Q}{I_iQ} \iso \frac{P}{I_iP} \quad {\rm be~isomorphism}~\ni\quad \det(\alpha_i)=\iota \otimes \frac{A}{I_i}.
$$
Let $\alpha:\frac{Q}{I_1I_2Q} \iso \frac{P}{I_1I_2P}$ be the isomorphism obtained by combining $\alpha_1$ and $\alpha_2$. Now, it follows 
that $\chi(\iota)(x)+\chi(\iota)(y)$ is obtained by combining $(I_1, \alpha_1\omega_1))$ and $(I_2, \alpha_2\omega_2))$, which is same as 
$\chi(\iota)(x+y)$.

\pic $\eop$


\bP\label{coDim0Funct} 
Use the notations as in Definition \ref{defnatMap}. Then,
Let   $Q'$ be another another projective $A$ module with $rank(Q)=d$ and let $\zeta:\Lambda^dQ' \iso \Lambda^dQ$ be an isomorphism.
Then, we have
$\chi(\iota\zeta)=\chi(\zeta) \chi(\iota)$. 
\eP

\pf Obvious! $\eop$
\bT\label{ItISGrp}
Use the notations as in Definition \ref{defnatMap}. Then, the map 
$\chi(\iota):\pi_0\left({\mathcal LO}(P)\right) \lra \pi_0\left({\mathcal LO}(Q)\right)$
 is a monoid isomorphism. 
 In particular, taking $Q=\Lambda^dQ \oplus A^{d-1}$, it follows $\pi_0\left({\mathcal LO}(P)\right)$ is a group.
\eT
\pf It follows from Proposition \ref{coDim0Funct}, that $\chi(\iota)\chi(\iota^{-1})=Id$ and $\chi(\iota^{-1})\chi(\iota)=Id$, Therefore, $\chi(\iota)$ 
is an isomorphism. 
Recall $\pi_0\left({\mathcal LO}(\Lambda^dQ \oplus A^{d-1})\right)$ is a group \cite[Theorem 6.11]{MM1}, which settles the latter part.
$\eop$

\begin{remark}{\rm 
Use the notations as in Definition \ref{defnatMap}. Theorem \ref{ItISGrp} asserts, $\pi_0\left({\mathcal LO}(P)\right)$ is a group, even when $P$ does not 
have a free direct summand. However, it remains open, whether $\pi_0\left({\mathcal LO}(P)\right)$ would fail to be a group, for some 
  projective $A$-module $P$, with $2rank(P)\geq d+2$.
}
\end{remark}

\bC\label{propChiD}
Let $A$ be an essentially smooth ring over an infinite perfect field $k$, with $1/2\in k$, and $\dim A=d\geq 3$. Let $Q$ be a projective $A$-module,
with $rank(Q)=n$. Let $P=\Lambda^d Q \oplus A^{d-1}$. Recall that that Euler class groups $E(P)\cong E^d(A, \Lambda^d Q)$,
as defined in \cite{MM1} and \cite{BS1}, respectively. The Euler class $e(Q) \in E^d(A, \Lambda^d Q)$, was also defined in \cite{BS1}.
 We we have a commutative diagram of isomorphisms
$$
\diagram 
\pi_0\left({\mathcal LO}(Q) \right) \ar[dr]\ar[r]^{\chi\qquad} & \pi_0\left({\mathcal LO}(\Lambda^d Q \oplus A^{d-1}) \right)\ar[d]_{\wr}^{\psi}\\
& E^d(A, \Lambda^d Q)\\
\enddiagram
$$
where the horizontal map $\chi$ is as in (\ref{defnatMap}), and the vertical isomorphism was defined in \cite[\S 7]{MM1} is induced by the the vertical isomorphism. 
We also have
$$
\psi\chi(\varepsilon(Q))=e(Q)
$$
where 
 $\varepsilon(Q):=[(0,0)]\in \pi_0\left({\mathcal LO}(Q) \right)$ denote the Nori homotopy class of $Q$, as defined in \cite{MM1}.
\eC
\pf Let $f:Q \sur I$ be a surjective homomorphism, where $I$ is an ideal of height $d$. Let $\omega: Q \sur I/I^2$ is induced by $f$.
Then, by definition $\varepsilon(Q)=[(I, \omega)]$. Rest of the proof follows from definition of Euler class $e(Q)$. 
\pic $\eop$

\section{Natural maps to Chow groups}\label{ChowSec}

The Chow groups of codimension  $n$ cycles will be denoted by
 $CH^n(A)$  (see \cite[\S 1.3]{F}). For an 
ideal $J\subseteq A$, $cycle\left(J\right)$ would denote the cycle of the closed set $V(J)$ (see \cite[\S 1.5]{F}).

\bD\label{DefMap2Chwo}{\rm 
Assume $A$ is a Cohen Macaulay ring, with $\dim A=d$. Also assume $\dim A_{\m}=d$ for all $\m\in \max(A)$.
As before, $P$ is a projective $A$-module with $rank(P)=n$. Define a map
$$
\ell^0_P: {\mathcal LO}^{\geq n}(P) \lra CH^n(A) \qquad {\rm by }\qquad \ell_P^0(I, \omega)=cycle(I)\in CH^n(A)
$$
Now define  the map 
$$
\ell_P: \pi_0\left({\mathcal LO}^{\geq n}(P)\right) \lra CH^n(A) \qquad {\rm by }\qquad \ell_P[(I, \omega)]=cycle(I)\in CH^n(A).
$$
We prove next that $\ell_P$ is a well defined set theoretic map.
}
\eD

\bP\label{propEllgeq}
Use the notations, as in (\ref{DefMap2Chwo}). Then, 
$$
\ell_P: \pi_0\left({\mathcal LO}^{\geq n}(P)\right) \lra CH^n(A)
$$
 is a well defined set theoretic map.
\eP
\pf 
For $t\in A$, let $\tau_t: CH^n(A[T]) \lra CH^n(A)$ denote the map induced by substitution $T=t$. Let $\pi: \spec{A}\lra \spec{A[T]}$ denote the 
structure map induced by $A\lra A[T]$ and $\pi^*: CH^n(A) \lra CH^n(A[T])$ denote the pull back map. Then, for any $t\in A$,  we have
$\tau_t=(\pi^*)^{-1}$  \TCP{\cite[Cor 6.5, pp. 111]{F}}. 

Let $(I_0, \omega_0), (I_1, \omega_1){\mathcal LO}^{\geq n}(P)$ and assume $(I_0, \omega_0)\sim (I_1, \omega_1)$. 
We prove $cycle(I_0)=cycle(I_1) \in CH^n(A)$.
There is a homotopy $H(T)=(J, \omega) \in 
{\mathcal LO}(P[T])$ such that $H(0)=(I_0, \omega_0)$ and $H(1)=(I_1, \omega_1)$. By Moving Lemma Argument (see Remark \ref{moviongRem}), we can assume 
$height(J)\geq n$.
Now, we have
\begin{equation}\label{TwoFulton}
\left\{
\begin{array}{l}
\ell_P^0(I_0, \omega_0) 
=cycle(J(0))=\tau_0(cycle(J))= (\pi^*)^{-1}(cycle(J))\\
\ell_P^0(I_1, \omega_1) 
= cycle(J(1))=\tau_1(cycle(J))= (\pi^*)^{-1}(cycle(J))\\
\end{array}
\right.
\end{equation}
\pic $\eop$
. 
\bT\label{finalCHThm}
Suppose $A$ is a regular ring, containing a filed $k$, with $1/2\in k$, and $\dim A=d$. Let $P$ be a projective $A$-module with $rank(P)=n$.
Then,
\bE
\item By Proposition \ref{transIso}, $\pi_0\left({\mathcal LO}^{\geq n}(P)\right)\iso \pi_0\left({\mathcal LO}(P)\right)$ is a bijection. Therefore,
 $\ell_P$ defined in (\ref{DefMap2Chwo}), defines a set theoretic map 
$$
\ell_P: \pi_0\left({\mathcal LO}(P)\right) \lra CH^n(A)
$$
\item  
Further,
$$
\ell_P(\varepsilon(P))=C^n(P^*)\in CH^n(A)
$$
where $\varepsilon(P)=[(0, 0)]\in \pi_0\left({\mathcal LO}(P)\right)$ denotes the Nori homotopy class   of $P$, and
 $C^n(P^*)$ denotes the top Chern class of the dual $P^*$.
\item Assume $2n\geq d+2$. 
Then, $\ell_P$ is an additive map.
\eE
\eT
\pf 
For the last statement, recall under the hypotheses $\pi_0\left({\mathcal LO}(P)\right)$, has structure of an abelian monoid. Rest of the proofs are obvious. $\eop$
\end{document}